# A remark about Galerkin method

N. Temirgaliyev[1]

1. **Introduction.** Galerkin method [1] which were created in 1915 has a powerful potential for practical applications [2]: "Any problem for which it is possible to write equations can be solved with one of the varieties of the Galerkin method".

At the same time, according to the Banach-Steinhaus resonance theorem, for every operator equation and for every linear method of solving this equation, including Galerkin method, whose linear bounded approximations have an unbounded sequence of norms, there always exists a right-hand side of the equation for which the sequence of approximations along this Method does not converge to its solution.

In this article, this effect is strengthened in the sense that any linear equation $Lu = f$ in the case of any linear versions of the Galerkin method, has at least as many unsolved right-hand sides in the form of linear combinations

$$f = L\psi_1 + ... + L\psi_N + L\psi_{N+1} + ... + L\psi_{5N} + ... + L\psi_T,$$

as there are finite-dimensional linear subspaces with dimensionality as much than five times as the number of basis functions $\psi_1,...,\psi_N$.

**2. Statement of the problem and formulation of the main result.** The Galerkin method in the linear case is as follow.

Let $L$ - linear operator from $X$ to $Y$, where $X$ and $Y$ - Banach spaces of functions which defined on $\Omega_X$ and $\Omega_Y$ respectively.

Let consider the equation

$$Lu = f. \qquad (1)$$

The source in the method of Galerkin are two systems: $\{\psi_k\}_{k=1}^{\infty}$ - linearly independent complete sequence (also called as *basic*) of functions in $X$ (as a rule, forming the basis of the space $X$) and $\{l_\tau\}_{\tau=1}^{\infty}$ - the complete sequence of linearly continuous functional on $Y$.

For given $N \geq 1$ approximate solution of equation (1) is sought in the form

$$u_N = u_N(x; f) = \sum_{k=1}^{N} c_k \psi_k(x). \qquad (2)$$

Then

$$Lu_N = \sum_{k=1}^{N} c_k L\psi_k,$$

where unknown $c_1,...,c_N$ are found from the condition

[1] Institute of theoretical mathematics and scientific computations, L.N.Gumilyov Eurasian National University, Astana, Kazakhstan, ntmath10@mail.ru

$$l_\tau(Lu_N - f) = 0 \ (\tau = 1,...,N),$$

i.e closeness $Lu_N$ to $f$ "measured" in the form of vanishing on the "discrepancy" $Lu_N - f$ of all linear functional $l_\tau \ (\tau = 1,...,N)$:

$$c_k(f) = \frac{\begin{vmatrix} l_1(L(\psi_1)) ... \overset{k-\text{ый столбец}}{l_1(f)} ... l_1(L(\psi_N)) \\ ... \\ l_N(L(\psi_1)) ... l_N(f) ... l_N(L(\psi_N)) \end{vmatrix}}{\det(l_\tau(L(\psi_k)))_{\tau,k=1}^N} \ (k = 1,...,N), \qquad (3)$$

on condition $\det(l_\tau(L(\psi_k)))_{\tau,k=1}^N \neq 0$.

Note that all coefficients $c_k(f) \ (k = 1,...,N)$ are linear functionals, because $l_1,...,l_N$ are linear.

Thus, the required approximate solution (2) of equation $Lu = f$ has form $u_N(f) = \sum_{k=1}^{N} c_k(f)\psi_k$, where $c_k(f)$ - linear functional, which is actually composed Galerkin method.

According to the definition of the completeness of the set of linear functional $\{l_\tau\}_{\tau=1}^\infty$, if for some $\bar{u}$ from $X$ and for all $\tau = 1,2,...$ is performed equality of "orthogonality"

$$\langle L\bar{u} - f, l_\tau \rangle \equiv l_\tau(L\bar{u} - f) = 0,$$

then $L\bar{u} = f$, i.e $\bar{u}$ is the exact solution of the equation $Lu = f$. Therefore, the Galerkin method can be understood as the "heuristic rule" based on this property: if for enough large $N$ holds the equalities

$$l_\tau(Lu_N - f) = 0 \ (\tau = 1,...,N), \qquad (4)$$

then norm of error $\|Lu_N - f\|_Y$ is sufficiently small.

In this connection the following question arises: how are connected execution of equation (4) for the given $N$ and the value (in terms of "small") of the $\|Lu_N - f\|_Y$?

A positive answer to this question must be as follows: for each $\varepsilon > 0$ there exists the positive integer $N = N_\varepsilon$, linearly independent on $X$ functions $\{\psi_k(x)\}_{k=1}^N (x \in \Omega_X)$ and functionals $\{l_\tau\}_{\tau=1}^N$ over $Y$, such that $\det(l_\tau(L(\psi_k)))_{\tau,k=1}^N \neq 0$ (this property is called the "correctness of the method") and for any function $f$ from Y with norm $\|f\|_Y \leq 1$ and built for her by the rules (2)-(3) functions $u_N(f) = \sum_{k=1}^{N} c_k(f)\psi_k$ satisfies equality

$$l_\tau(L(u_N(f))-f)=0 \ (\tau=1,\ldots,N)$$

and inequality

$$\|L(u_N(f))-f\|<\varepsilon.$$

The purpose of this article consists in conclusion that such an expectation is not justified.

First, we make some clarification in terms of the theorem stated below. According to the linear operator $L$ appropriate system $\{\psi_k\}_{k=1}^N$ of the linear independent functions from $X$ will assume so, that there is a linear independent in $Y$ system of function $\{\vartheta_k\}_{k=1}^N$ so, that and $\vartheta_k$ is continuous on the closure $\Omega_Y$ for all $k=1,\ldots,N$.

By $Y$ we understand space $L^2(\Omega_Y)$ summed in a square Lebesgue measurable functions defined on a set of Lebesgue $\Omega_X$ with a finite measure $|\Omega_Y|$.

Under these assumptions the following theorem is true.

**The main Theorem.** *Let a Banach space of real-valued functions $X$ and $Y=L^2(\Omega_Y)$ with set of definition $\Omega_X$ and $\Omega_Y$ respectively, and let $L$- linear operator from $X$ to $Y$.*

*Let be a positive integer $N\geq 2$ also given and linear independent system of functions $\{\psi_k\}_{k=1}^N$ from $X$, linear independent system of functional $\{c_k\}_{k=1}^N$ and $\{l_\tau\}_{\tau=1}^N$ over $Y$. Then will be find function $f_N$ from $L^2(\Omega_Y)$ with norm $\|f_N\|_{L^2(\Omega_Y)}=1$ so, that for function*

$$u_N(f_N)\equiv u_N(x;f_N)=\sum_{k=1}^N c_k(f_N)\psi_k(x) \qquad (5)$$

*is performed the equality*

$$l_\tau(L(u_N(f_N)))=l_\tau(f_N) \ (\tau=1,\ldots,N), \qquad (6)$$

*but*

$$\|f_N-L(u_N(f_N))\|_{L^2(\Omega_Y)}=1. \qquad (7)$$

**3. Remarks and conclusions.** $1^0$. For any linear operator $L$, all possible linear versions of the Galerkin method (for a given $N\geq 1$) are enclosed in a pair of linear aggregates of approximation

$$u_N(x;f_N)=\sum_{k=1}^N c_k(f_N)\psi_k(x), \quad x\in\Omega$$

and linear functionals $l_1,\ldots,l_N$ over $Y$.

We note in particular that, in the main theorem, the choice of linear functionals $\{c_k\}$ does not depend on systems $\{\psi_k\}_{k=1}^\infty$ and $\{l_\tau\}_{\tau=1}^N$, while in the Galerkin method they are connected by equalities (3).

$2^0$. The imposed conditions of linear independence $\{\vartheta_k\}_{k=1}^{N}$, should suppose, it is not particularly burdensome, since this property is across the boundary or other conditions provided by the uniqueness of solutions of the equation $Lu = f$. Conditions of continuity of $\vartheta_k$ improve the quality of the approximate solution.

$3^0$. Function $f_N$ which appropriate equation $Lu = f_N$ according to the Galerkin method is not solved approximately is a linear combination of the original functions $\vartheta_1(y) = L\psi_1, \ldots, \vartheta_N(y) = L\psi_N$ and arbitrary functions $\vartheta_{N+1}(y) = Lb_{N+1}, \ldots, \vartheta_T(y) = Lb_T$ in an arbitrary finite amount $T \geq 5N$ satisfying only the requirement of linear independence of the system $\{\psi_k\}_{k=1}^{T}$ in $X$ and $\{\upsilon_k\}_{k=1}^{T}$ in $Y$.

Let's make some clarification on the function $f_N$.

Firstly, $f_N$ located outside the linear span of the basic system $\{\vartheta_k\}_{k=1}^{N}$, which, apparently, ensures the existence of cases of effectiveness of the Galerkin method.

Secondly, if for construction $\bar{f}$ would require to involve linearly independent functions $\vartheta_k = L\psi_k$ in an infinite number, then

$$f_N \overset{Y}{=} \sum_{k=1}^{\infty} c_k \vartheta_k,$$

the function $f_N$ does not belong to any finite-dimensional subspace with basis from $\{\vartheta_k\}_{k=1}^{\infty}$, which would be natural and explained the Fundamental Theorem.

On this background, first, affiliation $f_N$ of the linear span of the finite linearly independent system $\{\vartheta_k\}_{k=1}^{T} (T \geq 5N)$ forms an unexpected effect.

And, finally, various functions $f_N$ as much as for all $T$, $5N \leq T < \infty$ for given $\vartheta_1(y) = L\psi_1, \ldots, \vartheta_N(y) = L\psi_N$ is possible to construct functions $\vartheta_{N+1}(y) = L\psi_{N+1}, \ldots, \vartheta_T(y) = L\psi_T$, collectively forming a linearly independent system $\psi \equiv \{\psi_k\}_{k=1}^{T}$ in $X$ and $\{\vartheta_k\}_{k=1}^{T}$ in $Y$.

$4^0$. The exact negation of "heuristic rules" formulated above for the method of Galerkin is: there is $\varepsilon_0 > 0$ that for any $N$, $\{\psi_k\}_{k=1}^{N}$ and $\{l_\tau\}_{\tau=1}^{N}$ such that $\det(l_\tau(L(\psi_k)))_{\tau,k=1}^{N} \neq 0$, will find $f_N$ with norm $\|f_N\|_Y \leq 1$ such that for function

$$\bar{u}_N \equiv \bar{u}_N\left(x, f_N; \{\psi_k\}_{k=1}^{N}, \{l_\tau\}_{\tau=1}^{N}\right) = \sum_{k=1}^{N} \frac{\begin{vmatrix} l_1(L(\psi_1)) & \ldots & \overset{k-\text{ый столбец}}{l_1(f_N)} & \ldots & l_1(L(\psi_N)) \\ & & \ldots & & \\ l_N(L(\psi_1)) & \ldots & l_N(f_N) & \ldots & l_N(L(\psi_N)) \end{vmatrix}}{\det(l_\tau(L(\psi_k)))_{\tau,k=1}^{N}} \psi_k(x) \quad (8)$$

is performed

$$l_\tau(L\bar{u}_N) = l_\tau(f_N) \quad (\tau = 1,...,N)$$

but

$$\|L\bar{u}_N - f_N\|_Y \geq \varepsilon_0. \tag{9}$$

The Fundamental Theorem states much more.

Firstly, (5) is replaced by $u_N(x) = \sum_{k=1}^{N} c_k(f_N)\psi_k(x)$ with an arbitrary set of linear functional $\{c_k\}_{k=1}^{N}$.

Secondly, even if would $\det\left[l_\tau\left(L(\psi_k)\right)\right]_{\tau,k=1}^{N} = 0$, the arbitrary choice of functional $\{c_k\}_{k=1}^{N}$ in (2) is still not ensure the effectiveness of the Galerkin method.

According to The Fundamental Theorem in (6) we can always take $\varepsilon_0 = 1$.

$5^0$. The Fundamental Theorem applies to the "*finite element method*" as contained in the "*Galerkin method*" when the area $\Omega_X$ is divided into a finite number of sub-areas that have the most common parts of the borders, and the basic functions are carriers of one or more of these sub-areas.

$6^0$. If in The Fundamental Theorem a linear operator $L$ is bounded with norm $\|L\|_{X \mapsto Y}$ and the equation $Lu = f_N$ has an exact solution $u_N^* \equiv u^*(f_N)$, then

$$\|u_N(f_N) - u_N^*\|_X \geq \frac{1}{\|L\|_{X \mapsto Y}}.$$

$7^0$. If $X = Y$ the Fundamental Theorem is applicable to solving both problems on the proper functions $u_\lambda$ and proper numbers $\lambda$ of the linear operator $Lu$ with Galerkin method, which is enough to proceed to the line operator $A_\lambda u = Lu - \lambda u$ and to equality $A_\lambda u = f$ when $f = 0$.

$8^0$. Galerkin method fits into the schemes of Computational (numerical) diameter (C(N)D) (see [3]). That is, for a given linear operator $L$ by $D_N$ we denote the conditions of main theorem.

Then in the notation adopted in [3], we have

$$\delta_N(0;L;Y_L;D_N)_{L^2(\Omega_Y)} \equiv \inf_{D_N} \sup_{f \in Y_L, \|f\|_Y \leq 1} \|f - L(u_N(f))\|_{L^2(\Omega_Y)} = \inf_{D_N} \sup_{f \in Y_L, \|f\|_Y \leq 1} \left\|f - \sum_{k=1}^{N} c_k(f)\vartheta_k\right\|_{L^2(\Omega_Y)}.$$

**The consequence of the Main Theorem.** *Rightly inequality*

$$\delta_N(0;L;Y_L;D_N)_{L^2(\Omega_Y)} \geq 1. \tag{10}$$

By (10) further advancement according to the scheme C(N)D in research portion C(N)D-2 and C(N)D-3 problems have no meaning.

By this C(N)D-approach warns of the presence of equation $Lu = f$ in a set of $\{f : f \in Y_L, \|f\|_Y \leq 1\}$ that is poorly solved by Galerkin method.